# Hypersurface family with a common geodesic

Ergin BAYRAM [a] and Emin KASAP [b]

## Abstract

In this paper, we study the problem of finding a hypersurface family from a given spatial geodesic curve in $\mathbb{R}^4$. We obtain the parametric representation for a hypersurface family whose members have the same curve as a given geodesic curve. Using the Frenet frame of the given geodesic curve, we present the hypersurface as a linear combination of this frame and analyze the necessary and sufficient condition for that curve to be geodesic. We illustrate this method by presenting some examples.

## 1. Introduction

Geodesic is a well-known notion in differential geometry. A geodesic on a surface can be defined in many equivalent ways. Geometrically, the shortest path joining any two points of a surface is a geodesic. Geodesics are curves in surfaces that play a role analogous that of straight lines in the plane. A straight line doesn't bend to left or right as we travel along it [6].

In recent years, there have been various researches on geodesics. Kumar et al., [20] presented a study on geodesic curves computed directly on NURBS surfaces and discrete geodesics computed on the equivalent tessellated surfaces. Wang et al., [26] studied the problem of constructing a family of surfaces from a given spatial geodesic curve and derived a parametric representation for a surface pencil whose members share the same geodesic curve as an isoparametric curve. A practical method was presented by Sanchez and Dorado, [21] to construct polynomial surfaces from a polynomial geodesic or a family of geodesics, by prescribing tangent ribbons. Sprynski et al., [22] dealt with reconstruction of numerical or real surfaces based on the knowledge of some geodesic curves on the surface. Paluszny, [19] considered patches that contain any given 3D polynomial curve as a pregeodesic (i.e. geodesic up to reparametrization). Given two pairs of regular space curves $\mathbf{r}_1(u)$, $\mathbf{r}_3(u)$ and $\mathbf{r}_2(v)$, $\mathbf{r}_4(v)$ that define a curvilinear rectangle, Farouki et al., [10] handled the problem of constructing a $C^2$ surface patch $\mathbf{R}(u,v)$ for which these four boundary curves correspond to geodesics of the surface. Farouki et al., [11] considered the problem of constructing polynomial or rational tensor-product Bézier patches bounded by given four polynomial or rational Bézier curves defining a curvilinear rectangle, such that they are geodesics of the constructed surface.

On the other hand, Wang et al., [26] tackled the problem of finding surfaces passing through a given geodesic. In 2011, given curve was changed to a line of curvature and Li et al., [18] constructed a surface family from a given line of curvature. Bayram et al., [5] gave the necessary and sufficient conditions for a given curve to be an asymptotic on a surface.

However, while differential geometry of a parametric surface in $\mathbb{R}^3$ can be found in textbooks such as in Struik [24], Willmore [28], Stoker [23], do Carmo [7], differential geometry of a parametric surface in $\mathbb{R}^n$ can be found in textbook such as in the contemporary literature on Geometric Modeling [9, 16]. Also, there is little literature on differential geometry of parametric surface family in $\mathbb{R}^3$ [2, 8, 17, 26], but not in $\mathbb{R}^4$. Besides, there is an ascending interest on fourth dimension [1, 2, 8].

Furthermore, various visualization techniques about objects in Euclidean n-space $(n \geq 4)$ are presented [3, 4, 14]. The fundamental step to visualize a 4D object is projecting first in to the 3-space and then into the plane. In many real world applications the problem of visualizing three-dimensional data, commonly referred to as scalar fields arouses. The graph of a function $\mathbf{f}(x,y,z): U \subset \mathbb{R}^3 \to \mathbb{R}$, U is open, is a special type of parametric hypersurface with the parametrization $(x, y, z, \mathbf{f}(x,y,z))$ in 4-space. There is an existing method for rendering such a 3-surface based on known methods for visualizing functions of two variables [13].

In this paper, we consider the four dimensional analogue of the problem of constructing a parametric representation of a surface family from a given spatial geodesic in Wang et al. [26], who derived the necessary and sufficient conditions on the marching-scale functions for which the curve C is an isogeodesic, i.e., both a geodesic and a parameter curve, on a given surface. We express the hypersurface pencil parametrically with the help of the Frenet frame $\{\mathbf{T}, \mathbf{N}, \mathbf{B}_1, \mathbf{B}_2\}$ of the given curve. We find the necessary and sufficient constraints on the marching-scale functions, namely, coefficients of Frenet vectors, so that both the geodesic and parametric requirements met. Finally, as an application of our method one example for each type of marching-scale functions is given.

## 2. Preliminaries

Let us first introduce some notations and definitions. Bold letters such as **a**, **R** will be used for vectors and vector functions. We assume that they are smooth enough so that all the (partial) derivatives given in the paper are meaningful. Let $\boldsymbol{\alpha}: \mathrm{I} \subset \mathbb{R} \to \mathbb{R}^4$ be an arc-length curve. If $\{\mathbf{T}, \mathbf{N}, \mathbf{B}_1, \mathbf{B}_2\}$ is the moving Frenet frame along $\boldsymbol{\alpha}$, then the Frenet formulas are given by

(1)
$$\begin{cases} \mathbf{T}' = k_1 \mathbf{N}, \\ \mathbf{N}' = -k_1 \mathbf{T} + k_2 \mathbf{B}_1, \\ \mathbf{B}_1' = -k_2 \mathbf{N} + k_3 \mathbf{B}_2, \\ \mathbf{B}_2' = -k_3 \mathbf{B}_1, \end{cases}$$

where $\mathbf{T}, \mathbf{N}, \mathbf{B}_1$ and $\mathbf{B}_2$ denote the tangent, principal normal, first binormal and second binormal vector fields, respectively, $k_i\ (i = 1,2,3)$ the i-th curvature functions of the curve $\boldsymbol{\alpha}$ [14].

From elementary differential geometry we have

(2)
$$\begin{cases} \boldsymbol{\alpha}'(s) = \mathbf{T}(s), \\ \boldsymbol{\alpha}''(s) = k_1(s)\mathbf{N}(s), \\ k_1(s) = \|\boldsymbol{\alpha}''(s)\|. \end{cases}$$

By using Frenet formulas one can obtain the followings

(3)
$$\begin{cases} \boldsymbol{\alpha}'''(s) = -k_1^2 \mathbf{T}(s) + k_1' \mathbf{N}(s) + k_1 k_2 \mathbf{B}_1(s), \\ \boldsymbol{\alpha}^{(iv)} = -3k_1 k_1' \mathbf{T}(s) + \left(-k_1^3 + k_1'' - k_1 k_2^2\right)\mathbf{N}(s) + \left(2k_1' k_2 + k_1 k_2'\right)\mathbf{B}_1(s) + k_1 k_2 k_3 \mathbf{B}_2(s). \end{cases}$$

The unit vectors $\mathbf{B}_2$ and $\mathbf{B}_1$ are given by

(4)
$$\begin{cases} \mathbf{B}_2(s) = \dfrac{\boldsymbol{\alpha}'(s) \otimes \boldsymbol{\alpha}''(s) \otimes \boldsymbol{\alpha}'''(s)}{\|\boldsymbol{\alpha}'(s) \otimes \boldsymbol{\alpha}''(s) \otimes \boldsymbol{\alpha}'''(s)\|}, \\ \mathbf{B}_1(s) = \mathbf{B}_2(s) \otimes \mathbf{T}(s) \otimes \mathbf{N}(s), \end{cases}$$

where $\otimes$ is the vector product of vectors in $\mathbb{R}^4$.

Since the vectors $\mathbf{T}$, $\mathbf{N}$, $\mathbf{B}_1$, $\mathbf{B}_2$ are orthonormal, the second curvature $k_2$ and the third curvature $k_3$ can be obtained from (3) as

$$(5)\qquad \begin{cases} k_2(s) = \dfrac{\mathbf{B}_1(s)\bullet\boldsymbol{\alpha}'''(s)}{k_1(s)}, \\ k_3(s) = \dfrac{\mathbf{B}_2(s)\bullet\boldsymbol{\alpha}^{(iv)}(s)}{k_1(s)k_2(s)}, \end{cases}$$

where ‘ • ’ denotes the standard inner product.

Let $\{\mathbf{e}_1, \mathbf{e}_2, \mathbf{e}_3, \mathbf{e}_4\}$ be the standard basis for four-dimensional Euclidean space $\mathbb{R}^4$. The vector product of the vectors $\mathbf{u} = \sum_{i=1}^{4} u_i\mathbf{e}_i, \mathbf{v} = \sum_{i=1}^{4} v_i\mathbf{e}_i, \mathbf{w} = \sum_{i=1}^{4} w_i\mathbf{e}_i$ is defined by

$$\mathbf{u}\otimes\mathbf{v}\otimes\mathbf{w} = \begin{vmatrix} \mathbf{e}_1 & \mathbf{e}_2 & \mathbf{e}_3 & \mathbf{e}_4 \\ u_1 & u_2 & u_3 & u_4 \\ v_1 & v_2 & v_3 & v_4 \\ w_1 & w_2 & w_3 & w_4 \end{vmatrix}$$

[15, 27].

If **u**, **v** and **w** are linearly independent then $\mathbf{u}\otimes\mathbf{v}\otimes\mathbf{w}$ is orthogonal to each of these vectors.

## 3. Hypersurface family with a common geodesic

A curve $\mathbf{r}(s)$ on a hypersurface $\mathbf{P} = \mathbf{P}(s,t,q) \subset \mathbb{R}^4$ is called an isoparametric curve if it is a parameter curve, that is, there exists a pair of parameters $t_0$ and $q_0$ such that $\mathbf{r}(s) = \mathbf{P}(s,t_0,q_0)$. Given a parametric curve $\mathbf{r}(s)$, it is called an isogeodesic of a hypersurface **P** if it is both a geodesic and an isoparametric curve on **P**.

Let $C:\mathbf{r} = \mathbf{r}\,(s),\ L_1 \le s \le L_2$, be a $C^3$ curve, where $s$ is the arc-length. To have a well-defined principal normal, assume that $\mathbf{r}''(s) \ne 0,\ L_1 \le s \le L_2$.

Let $\mathbf{T}(s), \mathbf{N}(s), \mathbf{B}_1(s), \mathbf{B}_2(s)$ be the tangent, principal normal, first binormal, second binormal, respectively; and let $k_1(s), k_2(s)$ and $k_3(s)$ be the first, second and the third curvature, respectively. Since $\{\mathbf{T}(s), \mathbf{N}(s), \mathbf{B}_1(s), \mathbf{B}_2(s)\}$ is an orthogonal coordinate frame on $\mathbf{r}(s)$ the parametric hypersurface $\mathbf{P}(s,t,q):[L_1,L_2]\text{x}[T_1,T_2]\text{x}[Q_1,Q_2] \to \mathbb{R}^4$ passing through $\mathbf{r}(s)$ can be defined as follows:

(6)
$$\mathbf{P}(s,t,q)=\mathbf{r}(s)+\big(\mathbf{u}(s,t,q),\mathbf{v}(s,t,q),\mathbf{w}(s,t,q),\mathbf{x}(s,t,q)\big)\begin{pmatrix}\mathbf{T}(s)\\ \mathbf{N}(s)\\ \mathbf{B}_1(s)\\ \mathbf{B}_2(s)\end{pmatrix}$$
$$L_1\le s\le L_2, T_1\le t\le T_2, Q_1\le q\le Q_2,$$

where $\mathbf{u}(s,t,q),\mathbf{v}(s,t,q),\mathbf{w}(s,t,q)$ and $\mathbf{x}(s,t,q)$ are all $C^1$ functions. These functions are called the *marching scale functions.*

We try to find out the necessary and sufficient conditions for a hypersurface $\mathbf{P}=\mathbf{P}(s,t,q)$ having the curve $C$ as an isogeodesic.

First to satisfy the isoparametricity condition there should exist $t_0\in[T_1,T_2]$ and $q_0\in[Q_1,Q_2]$ such that $\mathbf{P}(s,t_0,q_0)=\mathbf{r}(s), L_1\le s\le L_2$, that is,

(7)
$$\begin{cases}\mathbf{u}(s,t_0,q_0)=\mathbf{v}(s,t_0,q_0)=\mathbf{w}(s,t_0,q_0)=\mathbf{x}(s,t_0,q_0)\equiv 0,\\ t_0\in[T_1,T_2],\ q_0\in[Q_1,Q_2],\ L_1\le s\le L_2.\end{cases}$$

Secondly, the curve $C$ is a geodesic on the hypersurface $\mathbf{P}(s,t,q)$ if and only if the principal normal $\mathbf{N}(s)$ of the curve and the normal $\hat{\mathbf{n}}(s,t_0,q_0)$ of the hypersurface $\mathbf{P}(s,t,q)$ are linearly dependent, that is, parallel along the curve *C,* [25]. The normal $\hat{\mathbf{n}}(s,t_0,q_0)$ of the hypersurface can be obtained by calculating the vector product of the partial derivatives and using the Frenet formula as follows

$$\begin{aligned}\frac{\partial\mathbf{P}(s,t,q)}{\partial s}&=\left(1+\frac{\partial\mathbf{u}(s,t,q)}{\partial s}-\mathbf{v}(s,t,q)\kappa_1(s)\right)\mathbf{T}(s)\\&+\left(\mathbf{u}(s,t,q)\kappa_1(s)+\frac{\partial\mathbf{v}(s,t,q)}{\partial s}-\mathbf{w}(s,t,q)\kappa_2(s)\right)\mathbf{N}(s)\\&+\left(\mathbf{v}(s,t,q)\kappa_2(s)+\frac{\partial\mathbf{w}(s,t,q)}{\partial s}-\mathbf{x}(s,t,q)\kappa_3(s)\right)\mathbf{B}_1(s)\\&+\left(\mathbf{w}(s,t,q)\kappa_3(s)+\frac{\partial\mathbf{x}(s,t,q)}{\partial s}\right)\mathbf{B}_2(s),\end{aligned}$$

$$\begin{aligned}\frac{\partial\mathbf{P}(s,t,q)}{\partial t}&=\frac{\partial\mathbf{u}(s,t,q)}{\partial t}\mathbf{T}(s)+\frac{\partial\mathbf{v}(s,t,q)}{\partial t}\mathbf{N}(s)\\&+\frac{\partial\mathbf{w}(s,t,q)}{\partial t}\mathbf{B}_1(s)+\frac{\partial\mathbf{x}(s,t,q)}{\partial t}\mathbf{B}_2(s),\end{aligned}$$

and

$$\frac{\partial \mathbf{P}(s,t,q)}{\partial q}=\frac{\partial \mathbf{u}(s,t,q)}{\partial q}\mathbf{T}(s)+\frac{\partial \mathbf{v}(s,t,q)}{\partial q}\mathbf{N}(s)$$
$$+\frac{\partial \mathbf{w}(s,t,q)}{\partial q}\mathbf{B}_1(s)+\frac{\partial \mathbf{x}(s,t,q)}{\partial q}\mathbf{B}_2(s).$$

**Remark:** Because,

$$\begin{cases}\mathbf{u}(s,t_0,q_0)=\mathbf{v}(s,t_0,q_0)=\mathbf{w}(s,t_0,q_0)=\mathbf{x}(s,t_0,q_0)\equiv 0,\\ t_0\in[T_1,T_2],\ q_0\in[Q_1,Q_2],\ L_1\le s\le L_2.\end{cases}$$

along the curve *C,* by the definition of partial differentiation we have

$$\begin{cases}\dfrac{\partial \mathbf{u}(s,t_0,q_0)}{\partial s}=\dfrac{\partial \mathbf{v}(s,t_0,q_0)}{\partial s}=\dfrac{\partial \mathbf{w}(s,t_0,q_0)}{\partial s}=\dfrac{\partial \mathbf{x}(s,t_0,q_0)}{\partial s}\equiv 0,\\ t_0\in[T_1,T_2],\ q_0\in[Q_1,Q_2],\ L_1\le s\le L_2.\end{cases}$$

By using (7) we have

$$\hat{\mathbf{n}}(s,t_0,q_0)=\frac{\partial \mathbf{P}(s,t_0,q_0)}{\partial s}\otimes\frac{\partial \mathbf{P}(s,t_0,q_0)}{\partial t}\otimes\frac{\partial \mathbf{P}(s,t_0,q_0)}{\partial q}$$
$$=\phi_1(s,t_0,q_0)\mathbf{T}(s)-\phi_2(s,t_0,q_0)\mathbf{N}(s)+\phi_3(s,t_0,q_0)\mathbf{B}_1(s)-\phi_4(s,t_0,q_0)\mathbf{B}_2(s),$$

where

$$\phi_1(s,t_0,q_0)=\begin{vmatrix}\dfrac{\partial \mathbf{v}(s,t_0,q_0)}{\partial s} & \dfrac{\partial \mathbf{w}(s,t_0,q_0)}{\partial s} & \dfrac{\partial \mathbf{x}(s,t_0,q_0)}{\partial s}\\ \dfrac{\partial \mathbf{v}(s,t_0,q_0)}{\partial t} & \dfrac{\partial \mathbf{w}(s,t_0,q_0)}{\partial t} & \dfrac{\partial \mathbf{x}(s,t_0,q_0)}{\partial t}\\ \dfrac{\partial \mathbf{v}(s,t_0,q_0)}{\partial q} & \dfrac{\partial \mathbf{w}(s,t_0,q_0)}{\partial q} & \dfrac{\partial \mathbf{x}(s,t_0,q_0)}{\partial q}\end{vmatrix}=0$$

$$\phi_2(s,t_0,q_0)=\begin{vmatrix} 1+\dfrac{\partial \mathbf{u}(s,t_0,q_0)}{\partial s} & \dfrac{\partial \mathbf{w}(s,t_0,q_0)}{\partial s} & \dfrac{\partial \mathbf{x}(s,t_0,q_0)}{\partial s} \\ \dfrac{\partial \mathbf{u}(s,t_0,q_0)}{\partial t} & \dfrac{\partial \mathbf{w}(s,t_0,q_0)}{\partial t} & \dfrac{\partial \mathbf{x}(s,t_0,q_0)}{\partial t} \\ \dfrac{\partial \mathbf{u}(s,t_0,q_0)}{\partial q} & \dfrac{\partial \mathbf{w}(s,t_0,q_0)}{\partial q} & \dfrac{\partial \mathbf{x}(s,t_0,q_0)}{\partial q} \end{vmatrix}$$

$$=\begin{vmatrix} 1 & 0 & 0 \\ \dfrac{\partial \mathbf{u}(s,t_0,q_0)}{\partial t} & \dfrac{\partial \mathbf{w}(s,t_0,q_0)}{\partial t} & \dfrac{\partial \mathbf{x}(s,t_0,q_0)}{\partial t} \\ \dfrac{\partial \mathbf{u}(s,t_0,q_0)}{\partial q} & \dfrac{\partial \mathbf{w}(s,t_0,q_0)}{\partial q} & \dfrac{\partial \mathbf{x}(s,t_0,q_0)}{\partial q} \end{vmatrix}$$

$$=\frac{\partial \mathbf{w}(s,t_0,q_0)}{\partial t}\frac{\partial \mathbf{x}(s,t_0,q_0)}{\partial q}-\frac{\partial \mathbf{w}(s,t_0,q_0)}{\partial q}\frac{\partial \mathbf{x}(s,t_0,q_0)}{\partial t},$$

$$\phi_3(s,t_0,q_0)=\begin{vmatrix} 1+\dfrac{\partial \mathbf{u}(s,t_0,q_0)}{\partial s} & \dfrac{\partial \mathbf{v}(s,t_0,q_0)}{\partial s} & \dfrac{\partial \mathbf{x}(s,t_0,q_0)}{\partial s} \\ \dfrac{\partial \mathbf{u}(s,t_0,q_0)}{\partial t} & \dfrac{\partial \mathbf{v}(s,t_0,q_0)}{\partial t} & \dfrac{\partial \mathbf{x}(s,t_0,q_0)}{\partial t} \\ \dfrac{\partial \mathbf{u}(s,t_0,q_0)}{\partial q} & \dfrac{\partial \mathbf{v}(s,t_0,q_0)}{\partial q} & \dfrac{\partial \mathbf{x}(s,t_0,q_0)}{\partial q} \end{vmatrix}$$

$$=\begin{vmatrix} 1 & 0 & 0 \\ \dfrac{\partial \mathbf{u}(s,t_0,q_0)}{\partial t} & \dfrac{\partial \mathbf{v}(s,t_0,q_0)}{\partial t} & \dfrac{\partial \mathbf{x}(s,t_0,q_0)}{\partial t} \\ \dfrac{\partial \mathbf{u}(s,t_0,q_0)}{\partial q} & \dfrac{\partial \mathbf{v}(s,t_0,q_0)}{\partial q} & \dfrac{\partial \mathbf{x}(s,t_0,q_0)}{\partial q} \end{vmatrix}$$

$$=\frac{\partial \mathbf{v}(s,t_0,q_0)}{\partial t}\frac{\partial \mathbf{x}(s,t_0,q_0)}{\partial q}-\frac{\partial \mathbf{v}(s,t_0,q_0)}{\partial q}\frac{\partial \mathbf{x}(s,t_0,q_0)}{\partial t},$$

$$\phi_4(s,t_0,q_0)=\begin{vmatrix} 1+\dfrac{\partial \mathbf{u}(s,t_0,q_0)}{\partial s} & \dfrac{\partial \mathbf{v}(s,t_0,q_0)}{\partial s} & \dfrac{\partial \mathbf{w}(s,t_0,q_0)}{\partial s} \\ \dfrac{\partial \mathbf{u}(s,t_0,q_0)}{\partial t} & \dfrac{\partial \mathbf{v}(s,t_0,q_0)}{\partial t} & \dfrac{\partial \mathbf{w}(s,t_0,q_0)}{\partial t} \\ \dfrac{\partial \mathbf{u}(s,t_0,q_0)}{\partial q} & \dfrac{\partial \mathbf{v}(s,t_0,q_0)}{\partial q} & \dfrac{\partial \mathbf{w}(s,t_0,q_0)}{\partial q} \end{vmatrix}$$

$$=\begin{vmatrix} 1 & 0 & 0 \\ \dfrac{\partial \mathbf{u}(s,t_0,q_0)}{\partial t} & \dfrac{\partial \mathbf{v}(s,t_0,q_0)}{\partial t} & \dfrac{\partial \mathbf{w}(s,t_0,q_0)}{\partial t} \\ \dfrac{\partial \mathbf{u}(s,t_0,q_0)}{\partial q} & \dfrac{\partial \mathbf{v}(s,t_0,q_0)}{\partial q} & \dfrac{\partial \mathbf{w}(s,t_0,q_0)}{\partial q} \end{vmatrix}$$

$$=\frac{\partial \mathbf{v}(s,t_0,q_0)}{\partial t}\frac{\partial \mathbf{w}(s,t_0,q_0)}{\partial q}-\frac{\partial \mathbf{v}(s,t_0,q_0)}{\partial q}\frac{\partial \mathbf{w}(s,t_0,q_0)}{\partial t}.$$

So, $\hat{\mathbf{n}}(s,t_0,q_0)\,||\,\mathbf{N}(s)$ if and only if

$$
\text{(8)} \qquad \begin{array}{l} \phi_3(s,t_0,q_0) = \phi_4(s,t_0,q_0) \equiv 0,\ \phi_2(s,t_0,q_0) \neq 0, \\ t_0 \in [T_1,T_2],\ q_0 \in [Q_1,Q_2],\ L_1 \leq s \leq L_2. \end{array}
$$

Thus, any hypersurface defined by (6) has the curve $C$ as an isogeodesic if and only if

$$
\text{(9)} \qquad \begin{cases} \mathbf{u}(s,t_0,q_0) = \mathbf{v}(s,t_0,q_0) = \mathbf{w}(s,t_0,q_0) = \mathbf{x}(s,t_0,q_0) \equiv 0, \\ \phi_3(s,t_0,q_0) = \phi_4(s,t_0,q_0) \equiv 0,\ \phi_2(s,t_0,q_0) \neq 0, \\ t_0 \in [T_1,T_2],\ q_0 \in [Q_1,Q_2],\ L_1 \leq s \leq L_2, \end{cases}
$$

is satisfied. We call the set of hypersurfaces defined by (6) and satisfying (9) an *isogeodesic hypersurface family.*

To develop the method further and for simplification purposes, we analyze some types of marching-scale functions.

### 3.1. Marching-scale functions of type I

Let marching-scale functions be

$$
\begin{cases} \mathbf{u}(s,t,q) = \mathbf{l}(s)\mathbf{U}(t,q), \\ \mathbf{v}(s,t,q) = \mathbf{m}(s)\mathbf{V}(t,q), \\ \mathbf{w}(s,t,q) = \mathbf{n}(s)\mathbf{W}(t,q), \\ \mathbf{x}(s,t,q) = \mathbf{p}(s)\mathbf{X}(t,q), \end{cases} \qquad L_1 \leq s \leq L_2,\ T_1 \leq t \leq T_2, Q_1 \leq q \leq Q_2,
$$

where $\mathbf{l}(s),\mathbf{m}(s),\mathbf{n}(s),\mathbf{p}(s),\mathbf{U}(t,q),\mathbf{V}(t,q),\mathbf{W}(t,q),\mathbf{X}(t,q) \in C^1$ and $\mathbf{l}(s) \neq 0 \neq \mathbf{m}(s)$, $\mathbf{n}(s) \neq 0 \neq \mathbf{p}(s)$, $\forall s \in [L_1,L_2]$. By using (9) the necessary and sufficient condition for which the curve $C$ is an isogeodesic on the hypersurface $\mathbf{P}(s,t,q)$ can be given as

$$
\text{(10)} \qquad \begin{cases} \mathbf{U}(t_0,q_0) = \mathbf{V}(t_0,q_0) = \mathbf{W}(t_0,q_0) = \mathbf{X}(t_0,q_0) = 0, \\ \dfrac{\partial \mathbf{V}(t_0,q_0)}{\partial t}\dfrac{\partial \mathbf{X}(t_0,q_0)}{\partial q} - \dfrac{\partial \mathbf{V}(t_0,q_0)}{\partial q}\dfrac{\partial \mathbf{X}(t_0,q_0)}{\partial t} = 0, \\ \dfrac{\partial \mathbf{V}(t_0,q_0)}{\partial t}\dfrac{\partial \mathbf{W}(t_0,q_0)}{\partial q} - \dfrac{\partial \mathbf{V}(t_0,q_0)}{\partial q}\dfrac{\partial \mathbf{W}(t_0,q_0)}{\partial t} = 0, \\ \dfrac{\partial \mathbf{W}(t_0,q_0)}{\partial t}\dfrac{\partial \mathbf{X}(t_0,q_0)}{\partial q} - \dfrac{\partial \mathbf{W}(t_0,q_0)}{\partial q}\dfrac{\partial \mathbf{X}(t_0,q_0)}{\partial t} \neq 0, \end{cases}
$$

$$
t_0 \in [T_1,T_2],\ q_0 \in [Q_1,Q_2].
$$

With a closer investigation of (10), we should have $\frac{\partial \mathbf{V}(t_0,q_0)}{\partial t}=0$ and $\frac{\partial \mathbf{V}(t_0,q_0)}{\partial q}=0$.

So, (10) can be simplified to

$$(11)\qquad \begin{cases} \mathbf{U}(t_0,q_0)=\mathbf{V}(t_0,q_0)=\mathbf{W}(t_0,q_0)=\mathbf{X}(t_0,q_0)=0, \\ \dfrac{\partial \mathbf{V}(t_0,q_0)}{\partial t}=\dfrac{\partial \mathbf{V}(t_0,q_0)}{\partial q}=0, \\ \dfrac{\partial \mathbf{W}(t_0,q_0)}{\partial t}\dfrac{\partial \mathbf{X}(t_0,q_0)}{\partial q}-\dfrac{\partial \mathbf{W}(t_0,q_0)}{\partial q}\dfrac{\partial \mathbf{X}(t_0,q_0)}{\partial t}\neq 0, \end{cases}$$

$$t_0\in[T_1,T_2],\ q_0\in[Q_1,Q_2].$$

### 3.2. Marching-scale functions of type II

Let marching-scale functions be

$$\begin{cases} \mathbf{u}(s,t,q)=\mathbf{l}(s,t)\mathbf{U}(q), \\ \mathbf{v}(s,t,q)=\mathbf{m}(s,t)\mathbf{V}(q), \\ \mathbf{w}(s,t,q)=\mathbf{n}(s,t)\mathbf{W}(q), \\ \mathbf{x}(s,t,q)=\mathbf{p}(s,t)\mathbf{X}(q), \end{cases}\qquad L_1\le s\le L_2,\ T_1\le t\le T_2,\ Q_1\le q\le Q_2,$$

where $\mathbf{l}(s),\mathbf{m}(s),\mathbf{n}(s),\mathbf{p}(s),\mathbf{U}(t,q),\mathbf{V}(t,q),\mathbf{W}(t,q),\mathbf{X}(t,q)\in C^1$. Also let us choose

$\mathbf{V}(q_0)=\frac{d\mathbf{V}(q_0)}{dq}=\mathbf{U}(q_0)=\frac{d\mathbf{U}(q_0)}{dq}=0$. By using (9), the curve $C$ is an isogeodesic on the hypersurface $\mathbf{P}(s,t,q)$ if and only if the followings are satisfied

$$(12)\qquad \begin{cases} \mathbf{n}(s,t_0)\mathbf{W}(q_0)=\mathbf{p}(s,t_0)\mathbf{X}(q_0)\equiv 0, \\ \dfrac{\partial \mathbf{n}(s,t_0)}{\partial t}\mathbf{W}(q_0)\mathbf{p}(s,t_0)\dfrac{d\mathbf{X}(q_0)}{dq}-\mathbf{n}(s,t_0)\dfrac{d\mathbf{W}(q_0)}{dq}\dfrac{\partial \mathbf{p}(s,t_0)}{\partial t}\mathbf{X}(q_0)\neq 0, \end{cases}$$

$$t_0\in[T_1,T_2],q_0\in[Q_1,Q_2],L_1\le s\le L_2.$$

### 3.3. Marching-scale functions of type III

Let marching-scale functions be

$$\begin{cases} \mathbf{u}(s,t,q) = \mathbf{l}(s,q)\mathbf{U}(t), \\ \mathbf{v}(s,t,q) = \mathbf{m}(s,q)\mathbf{V}(t), \\ \mathbf{w}(s,t,q) = \mathbf{n}(s,q)\mathbf{W}(t), \\ \mathbf{x}(s,t,q) = \mathbf{p}(s,q)\mathbf{X}(t), \end{cases} \qquad L_1 \le s \le L_2,\ T_1 \le t \le T_2,\ Q_1 \le q \le Q_2,$$

where $\mathbf{l}(s,q),\ \mathbf{m}(s,q),\ \mathbf{n}(s,q),\ \mathbf{p}(s,q), \mathbf{U}(t),\ \mathbf{V}(t),\ \mathbf{W}(t),\ \mathbf{X}(t) \in C^1$. Also let $\mathbf{V}(t_0) = \frac{d\mathbf{V}(t_0)}{dt} = \mathbf{U}(t_0) = \frac{d\mathbf{U}(t_0)}{dt} = 0$. By using (9) we derive the necessary and sufficient condition for which the curve $C$ is an isogeodesic on the hypersurface $\mathrm{P}(s,t,q)$ as

(13)
$$\begin{cases} \mathbf{n}(s,q_0)\mathbf{W}(t_0) = \mathbf{p}(s,q_0)\mathbf{X}(t_0) \equiv 0, \\ \mathbf{n}(s,t_0)\frac{d\mathbf{W}(t_0)}{dt}\frac{\partial \mathbf{p}(s,q_0)}{\partial q}\mathbf{X}(t_0) - \frac{\partial \mathbf{n}(s,q_0)}{\partial q}\mathbf{W}(t_0)\mathbf{p}(s,q_0)\frac{d\mathbf{X}(t_0)}{dt} \ne 0, \end{cases}$$

$$t_0 \in [T_1,T_2],\ q_0 \in [Q_1,Q_2], L_1 \le s \le L_2.$$

## 4. Examples

**Example 1.** Let $\mathbf{r}(s) = \left(\frac{1}{2}\cos(s), \frac{1}{2}\sin(s), \frac{1}{2}s, \frac{\sqrt{2}}{2}s\right),\ 0 \le s \le 2\pi,$ be a curve parametrized by arc-length. For this curve,

$$\begin{aligned} &\mathbf{T}(s) = \mathbf{r}'(s) = \left(-\frac{1}{2}\sin(s), \frac{1}{2}\cos(s), \frac{1}{2}, \frac{\sqrt{2}}{2}\right), \\ &\mathbf{N}(s) = \left(-\cos(s), -\sin(s), 0, 0\right), \\ &\mathbf{B}_2(s) = \frac{r'(s)\otimes r''(s)\otimes r'''(s)}{\left\| r'(s)\otimes r''(s)\otimes r'''(s)\right\|} = \left(0, 0, \frac{\sqrt{6}}{3}, -\frac{\sqrt{3}}{3}\right), \\ &\mathbf{B}_1(s) = \mathbf{B}_2 \otimes \mathbf{T} \otimes \mathbf{N} = \left(-\frac{\sqrt{3}}{2}\sin(s), \frac{\sqrt{3}}{2}\cos(s), -\frac{\sqrt{3}}{6}, -\frac{\sqrt{6}}{6}\right). \end{aligned}$$

Let us choose the marching-scale functions of type I, where

$$\mathbf{l}(s) = \mathbf{m}(s) = \mathbf{n}(s) = \mathbf{p}(s) \equiv 1$$

and

$$\mathbf{U}(t,q) = (t-t_0)(q-q_0),\ \mathbf{V}(t,q) \equiv 0,\ \mathbf{W}(t,q) = t - t_0,\ \mathbf{X}(t,q) = q - q_0,$$

$t_0 \in [0,1], q_0 \in [0,1], 0 \le s \le 2\pi.$

So, we have

$$\begin{aligned}
&\mathbf{u}(s,t,q)=(t-t_0)(q-q_0),\\
&\mathbf{v}(s,t,q)\equiv 0,\\
&\mathbf{w}(s,t,q)=t-t_0,\\
&\mathbf{x}(s,t,q)=q-q_0.
\end{aligned}$$

The hypersurface

$$\begin{aligned}
\mathbf{P}(s,t,q)&=\mathbf{r}(s)+\mathbf{u}(s,t,q)\mathbf{T}(s)+\mathbf{v}(s,t,q)\mathbf{N}(s)+\mathbf{w}(s,t,q)\mathbf{B}_1(s)+\mathbf{x}(s,t,q)\mathbf{B}_2(s)\\
&=\left(\frac{1}{2}\cos(s)-\frac{1}{2}(t-t_0)(q-q_0)\sin(s)-\frac{\sqrt{3}}{2}(t-t_0)\sin(s),\right.\\
&\quad\frac{1}{2}\sin(s)+\frac{1}{2}(t-t_0)(q-q_0)\cos(s)+\frac{\sqrt{3}}{2}(t-t_0)\cos(s),\\
&\quad\frac{1}{2}s+\frac{1}{2}(t-t_0)(q-q_0)-\frac{\sqrt{3}}{6}(t-t_0)+\frac{\sqrt{6}}{3}(q-q_0),\\
&\quad\left.\frac{\sqrt{2}}{2}s+\frac{\sqrt{2}}{2}(t-t_0)(q-q_0)-\frac{\sqrt{6}}{6}(t-t_0)-\frac{\sqrt{3}}{3}(q-q_0)\right)
\end{aligned}$$

$0\le s\le 2\pi,\ 0\le t\le 1,\ 0\le q\le 1,\ t_0\in[0,1], q_0\in[0,1],$ is a member of the isogeodesic hypersurface family, since it satisfies (11).

By changing the parameters $t_0$ and $q_0$ we can adjust the position of the curve $\mathbf{r}(s)$ on the hypersurface. Let us choose $t_0=\frac{1}{2}$ and $q_0=0$. Now the curve $\mathbf{r}(s)$ is again an isogeodesic on the hypersurface $\mathbf{P}(s,t,q)$ and the equation of the hypersurface is

$$\begin{aligned}
\mathbf{P}(s,t,q)=&\left(\frac{1}{2}\cos(s)-\frac{1}{2}\left(t-\frac{1}{2}\right)\left(q+\sqrt{3}\right)\sin(s),\right.\\
&\frac{1}{2}\sin(s)+\frac{1}{2}\left(t-\frac{1}{2}\right)\left(q+\sqrt{3}\right)\cos(s),\\
&\frac{1}{2}s+\frac{1}{2}\left(t-\frac{1}{2}\right)q-\frac{\sqrt{3}}{6}\left(t-\frac{1}{2}\right)+\frac{\sqrt{6}}{3}q,\\
&\left.\frac{\sqrt{2}}{2}s+\frac{\sqrt{2}}{2}\left(t-\frac{1}{2}\right)q-\frac{\sqrt{6}}{6}\left(t-\frac{1}{2}\right)-\frac{\sqrt{3}}{3}q\right).
\end{aligned}$$

The projection of a hypersurface into 3-space generally yields a three-dimensional volume. If we fix each of the three parameters, one at a time, we obtain three distinct families of 2-spaces in 4-space. The projections of these 2-surfaces into 3-space are surfaces in 3-space. Thus, they can be displayed by 3D rendering methods.

So, if we (parallel) project the hypersurface $\mathbf{P}(s,t,q)$ into the $w=0$ subspace and fix $q=\frac{1}{8}$ we obtain the surface

$$\mathbf{P}_w\left(s,t,\frac{1}{8}\right)=\left(\frac{1}{2}\cos(s)-\frac{1+8\sqrt{3}}{16}\left(t-\frac{1}{2}\right)\sin(s),\right.$$
$$\frac{1}{2}\sin(s)+\frac{1+8\sqrt{3}}{8}\left(t-\frac{1}{2}\right)\cos(s),$$
$$\left.\frac{1}{2}s+\frac{1}{16}\left(t-\frac{1}{2}\right)-\frac{\sqrt{3}}{6}\left(t-\frac{1}{2}\right)+\frac{\sqrt{6}}{24}\right),\ 0\le s\le 2\pi,\ 0\le t\le 1$$

in 3-space illustrated in Fig. 1.

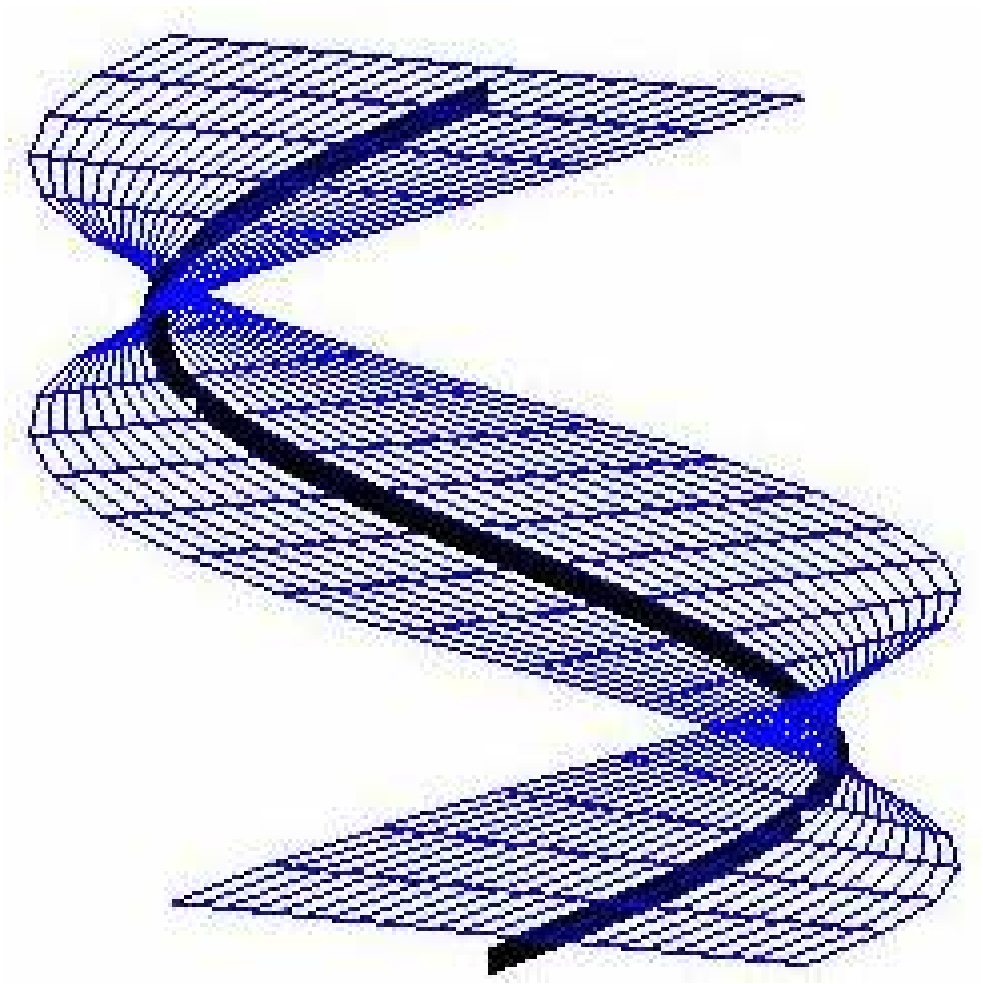

**Fig. 1**. Projection of a member of the hypersurface family with marching-scale functions of type I and its isogeodesic.

**Example 2.** Given the curve parameterized by arc-length

$\mathbf{r}(s)=\left(\frac{1}{2}\sin(s),\frac{1}{2}\cos(s),0,\frac{\sqrt{3}}{2}s\right),\ 0\le s\le 2\pi,$ it is easy to show that

$$
\begin{aligned}
&\mathbf{T}(s)=\mathbf{r}'(s)=\left(\frac{1}{2}\cos(s),-\frac{1}{2}\sin(s),0,\frac{\sqrt{3}}{2}\right),\\
&\mathbf{N}(s)=\left(-\sin(s),-\cos(s),0,0\right),\\
&\mathbf{B}_2(s)=\frac{\mathbf{r}'(s)\otimes\mathbf{r}''(s)\otimes\mathbf{r}'''(s)}{\left\|\mathbf{r}'(s)\otimes\mathbf{r}''(s)\otimes\mathbf{r}'''(s)\right\|}=\left(0,0,-1,0\right),\\
&\mathbf{B}_1(s)=\mathbf{B}_2\otimes\mathbf{T}\otimes\mathbf{N}=\left(\frac{\sqrt{3}}{2}\cos(s),-\frac{\sqrt{3}}{2}\sin(s),0,-\frac{1}{2}\right).
\end{aligned}
$$

Let us choose the marching-scale functions of type II, where

$$\mathbf{n}(s,t)=s+t+1,\ \mathbf{p}(s,t)=(s+1)(t-t_0),$$

and

$$\mathbf{U}(q)=\mathbf{V}(q)\equiv 0,\mathbf{W}(q)=q-q_0,\ \mathbf{X}(q)=1.$$

So, we get

$$
\begin{aligned}
&\mathbf{u}(s,t,q)\equiv 0,\\
&\mathbf{v}(s,t,q)\equiv 0,\\
&\mathbf{w}(s,t,q)=(s+t+1)(q-q_0),\\
&\mathbf{x}(s,t,q)=(s+1)(t-t_0).
\end{aligned}
$$

From (12) the hypersurface

$$
\begin{aligned}
\mathbf{P}(s,t,q)&=\mathbf{r}(s)+\mathbf{u}(s,t,q)\mathbf{T}(s)+\mathbf{v}(s,t,q)\mathbf{N}(s)+\mathbf{w}(s,t,q)\mathbf{B}_1(s)+\mathbf{x}(s,t,q)\mathbf{B}_2(s)\\
&=\Bigg(\frac{1}{2}\sin(s)+\frac{\sqrt{3}}{2}(s+t+1)(q-q_0)\cos(s),\\
&\quad\frac{1}{2}\cos(s)-\frac{\sqrt{3}}{2}(s+t+1)(q-q_0)\sin(s),\\
&\quad-(s+1)(t-t_0),\\
&\quad\frac{\sqrt{3}}{2}\mathrm{s}-\frac{1}{2}(s+t+1)(q-q_0)\Bigg)
\end{aligned}
$$

$0\le s\le 2\pi,\ 0\le t\le 1,\ 0\le q\le 1,$ is a member of the hypersurface family having the curve $\mathbf{r}(s)$ as an isogeodesic.

Setting $t_0=\frac{1}{2}$ and $q_0=0$ yields the hypersurface

$$\mathbf{P}(s,t,q)=\left(\frac{1}{2}\sin(s)+\frac{\sqrt{3}}{2}(s+t+1)q\cos(s),\right.$$
$$\frac{1}{2}\cos(s)-\frac{\sqrt{3}}{2}(s+t+1)q\sin(s),$$
$$-(s+1)\left(t-\frac{1}{2}\right),$$
$$\left.\frac{\sqrt{3}}{2}s-\frac{1}{2}(s+t+1)q\right).$$

By (parallel) projecting the hypersurface $\mathbf{P}(s,t,q)$ into the $\mathrm{w}=0$ subspace and fixing $q=\frac{1}{500}$ we get the surface

$$\mathbf{P}_w\left(s,t,\frac{1}{500}\right)=\left(\frac{1}{2}\sin(s)+\frac{\sqrt{3}}{1000}(s+t+1)\cos(s),\right.$$
$$\frac{1}{2}\cos(s)-\frac{\sqrt{3}}{1000}(s+t+1)\sin(s),$$
$$\left.-(s+1)\left(t-\frac{1}{2}\right)\right),$$

where, $0\le s\le 2\pi,\ 0\le t\le 1$ in 3-space demonstrated in Fig. 2.

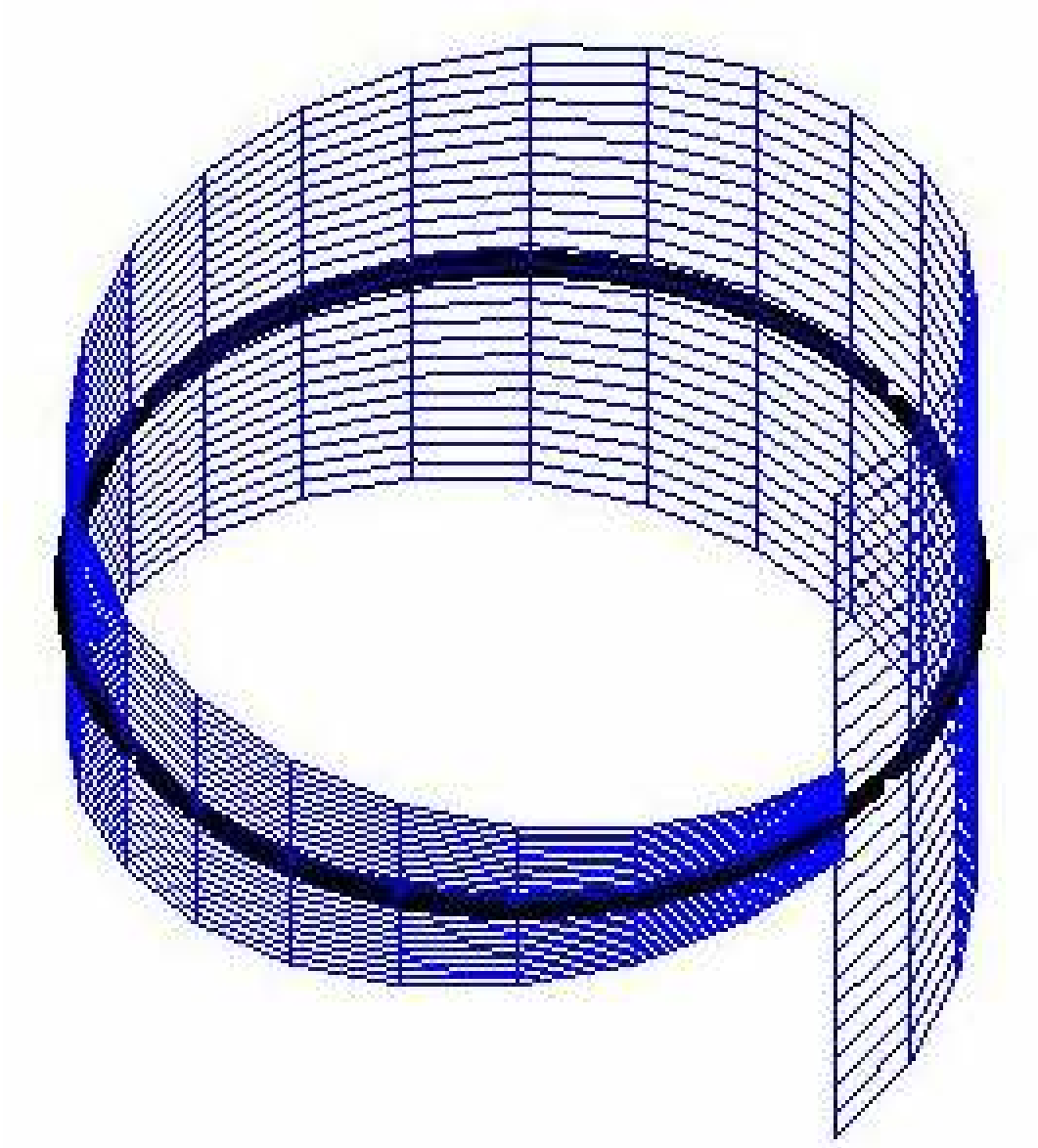

**Fig. 2.** Projection of a member of the hypersurface family with marching-scale functions of type II and its isogeodesic.

**Example 3.** Let $\mathbf{r}(s)=\left(\frac{1}{2}\sin(s),\frac{1}{2}\cos(s),0,\frac{\sqrt{3}}{2}s\right)$, $0\le s\le 2\pi$, be an arc-length curve. One can easily show that

$$\begin{aligned}
&\mathbf{T}(s)=\mathbf{r}'(s)=\left(\frac{1}{2}\cos(s),-\frac{1}{2}\sin(s),0,\frac{\sqrt{3}}{2}\right),\\
&\mathbf{N}(s)=\left(-\sin(s),-\cos(s),0,0\right),\\
&\mathbf{B}_2(s)=\frac{\mathbf{r}'(s)\otimes\mathbf{r}''(s)\otimes\mathbf{r}'''(s)}{\left\|\mathbf{r}'(s)\otimes\mathbf{r}''(s)\otimes\mathbf{r}'''(s)\right\|}=(0,0,-1,0),\\
&\mathbf{B}_1(s)=\mathbf{B}_2\otimes\mathbf{T}\otimes\mathbf{N}=\left(\frac{\sqrt{3}}{2}\cos(s),-\frac{\sqrt{3}}{2}\sin(s),0,-\frac{1}{2}\right),
\end{aligned}$$

for this curve.

If we choose the marching-scale functions of type III, where

$$\mathbf{n}(s,q)=\sin\left(s(q-q_0)\right),\ \mathbf{p}(s,q)=sq^2$$

and

$$\mathbf{U}(t)=\mathbf{V}(t)\equiv 0,\ \mathbf{W}(t)=1,\ \mathbf{X}(t)=t-t_0$$

then

$$\begin{aligned}
&\mathbf{u}(s,t,q)\equiv 0,\\
&\mathbf{v}(s,t,q)\equiv 0,\\
&\mathbf{w}(s,t,q)=\sin\left(s(q-q_0)\right),\\
&\mathbf{x}(s,t,q)=sq^2(t-t_0).
\end{aligned}$$

Thus, from (13) if we take $q_0\ne 0$ then the curve $\mathbf{r}(s)$ is an isogeodesic on the hypersurface

$$\begin{aligned}
\mathbf{P}(s,t,q)&=\mathbf{r}(s)+\underbrace{\mathbf{u}(s,t,q)}_{0}\mathbf{T}(s)+\underbrace{\mathbf{v}(s,t,q)}_{0}\mathbf{N}(s)+\mathbf{w}(s,t,q)\mathbf{B}_1(s)+\mathbf{x}(s,t,q)\mathbf{B}_2(s)\\
&=\Bigg(\frac{1}{2}\sin(s)+\frac{\sqrt{3}}{2}\cos(s)\sin\left(s(q-q_0)\right),\\
&\quad\frac{1}{2}\cos(s)-\frac{\sqrt{3}}{2}\sin(s)\sin\left(s(q-q_0)\right),\\
&\quad -sq^2(t-t_0),\\
&\quad\frac{\sqrt{3}}{2}s-\frac{1}{2}\sin\left(s(q-q_0)\right)\Bigg),
\end{aligned}$$

where $\pi\le s\le 3\pi,\ 0\le t\le 1,\ 0\le q\le 1.$

By taking $t_0 = 1$ and $q_0 = 1$ we have the following hypersurface:

$$\mathbf{P}(s,t,q) = \left( \frac{1}{2}\sin(s) + \frac{\sqrt{3}}{2}\cos(s)\sin(s(q-1)), \right.$$
$$\frac{1}{2}\cos(s) - \frac{\sqrt{3}}{2}\sin(s)\sin(s(q-1)),$$
$$-sq^2(t-1),$$
$$\left. \frac{\sqrt{3}}{2}\mathrm{s} - \frac{1}{2}\sin(s(q-1)) \right).$$

Hence, if we (parallel) project the hypersurface $\mathrm{P}(s,t,q)$ into the $z = 0$ subspace we get the surface

$$\mathbf{P}_z(s,q) = \left( \frac{1}{2}\sin(s) + \frac{\sqrt{3}}{2}\cos(s)\sin(s(q-1)), \right.$$
$$\frac{1}{2}\cos(s) - \frac{\sqrt{3}}{2}\sin(s)\sin(s(q-1)),$$
$$\left. \frac{\sqrt{3}}{2}s - \frac{1}{2}\sin(s(q-1)) \right),$$

where $\pi \le s \le 3\pi,\ 0 \le q \le 1,$

in 3-space shown in Fig. 3.

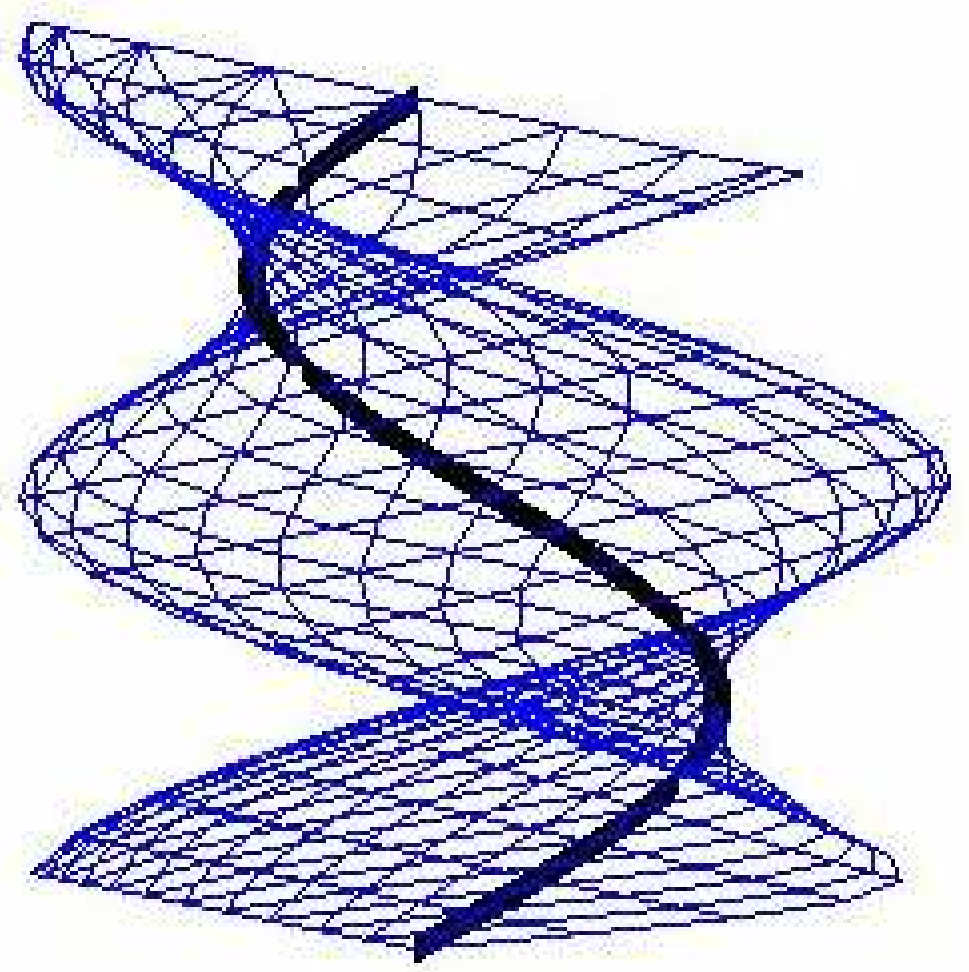

**Fig. 3.** Projection of a member of the hypersurface family with marching-scale functions of type III and its isogeodesic.

## 5. Conclusion

We have introduced a method for finding a hypersurface family passing through the same given geodesic as an isoparametric curve. The members of the hypersurface family are obtained by choosing suitable marching-scale functions. For a better analysis of the method we investigate three types of marching-scale functions. Also, by giving an example for each type the method is verified. Furthermore, with the help of the projecting methods a member of the family is visualized in 3-space with its isogeodesic.

However, there is more work waiting to study. For 3-space, one possible alternative is to consider the realm of implicit surfaces $\mathbf{F}(x,y,z,t)=0$ and try to find out the constraints for a given curve $\mathbf{r}(s)$ is an isogeodesic on $\mathbf{F}(x,y,z,t)=0$. Also, the analogue of the problem dealt in this paper may be considered for 2-surfaces in 4-space or another types of marching-scale functions may be investigated.

## 6. References

[1] Abdel-All N.H., Badr S.A., Soliman M.A., Hassan S.A., Intersection curves of hypersurfaces in $\mathbb{R}^4$, Comput. Aided Geomet. Design, 29, (2012), 99-108.

[2] Aléssio O., Differential geometry of intersection curves in $\mathbb{R}^4$ of three implicit surfaces, Comput. Aided Geom. Design, 26, (2009), 455–471.

[3] Banchoff T.F., Beyond the Third Dimension: geometry, computer graphics, and higher dimensions, W.H. Freeman & Co., New York, NY, USA, 1990.

[4] Banchoff T.F., Discovering the Fourth Dimension, Prime Computer, Inc., Natick, MA, 1987.

[5] Bayram E., Güler F., Kasap E., Parametric representation of a surface pencil with a common asymptotic curve, Comput. Aided Des., 44, (2012), 637-643.

[6] Bloch E.D., A First Course in Geometry, Birkhäuser, Boston, 1997.

[7] do Carmo M.P., Differential Geometry of Curves and Surfaces, Prentice Hall, Englewood Cliffs, NJ, 1976.

[8] Düldül M., On the intersection curve of three parametric hypersurfaces, Comput. Aided Geom. Design, 27, (2010), 118-127.

[9] Farin G., Curves and Surfaces for Computer Aided Geometric Design: a practical guide. Academic Press, Inc., San Diego, CA, 2002.

[10] Farouki R.T., Szafran N., Biard L., Existence conditions for Coons patches interpolating geodesic boundary curves, Comput. Aided Geom. Design, 26, (2009a), 599-614.

[11] Farouki R.T., Szafran N., Biard L., Construction of Bézier surface patches with Bézier curves as geodesic boundaries, Comput. Aided Geom. Design, 41, (2009b), 772-781.

[12] Gluck, H. Higher curvatures of curves in Euclidean space, Amer. Math. Monthly, 73, (1966), 699–704.

[13] Hamann B., Visualization and modeling contours of trivariate functions, Ph.D. thesis, Arizona State University, 1991.

[14] Hanson A.J., Heng P.A., Visualizing the fourth dimension using geometry and light, In: Proceedings of the 2nd Conference on Visualization '91, IEEE Computer Society Press, Los Alamitos, CA, USA, pp. 321–328, 1991.

[15] Hollasch S.R., Four-space visualization of 4D objects, Master thesis, Arizona State University, 1991.

[16] Hoschek J., Lasser D., Fundamentals of Computer Aided Geometric Design, A.K. Peters, Wellesley, MA, 1993.

[17] Kasap E., Akyildiz F.T., Orbay K., A generalization of surfaces family with common spatial geodesic, Appl. Math. Comput., 201, (2008), 781-789.

[18] Li C.Y., Wang R.H., Zhu C.G, Parametric representation of a surface pencil with a common line of curvature, Comput. Aided Des., 43, (2011), 1110-1117.

[19] Paluszny M., Cubic polynomial patches through geodesics, Comput. Aided Des., 40, (2008), 56-61.

[20] Ravi Kumar G.V.V., Srinivasan P., Holla V.D, Shastry K.G., Prakash B.G., Geodesic curve computations on surfaces, Comput. Aided Geom. Design 20, (2003), 119-133.

[21] Sánchez-Reyes J., Dorado R., Constrained design of polynomial surfaces from geodesic curves, Comput. Aided Des., 40, (2008), 49-55.

[22] Sprynski N., Szafran N., Lacolle B., Biard L., Surface reconstruction via geodesic interpolation, Comput. Aided Des., 40, (2008), 480-492.

[23] Stoker J.J., Differential geometry, Wiley, New York, 1969.

[24] Struik D.J., Lectures on Classical Differential Geometry, Addison–Wesley, Reading, MA, 1950.

[25] Thorpe J.A., Elementary Topics in Differential Geometry, Springer-Verlag, New York, Heidelberg-Berlin, 1979.

[26] Wang G.J., Tang K., Tai C.L., Parametric representation of a surface pencil with a common spatial geodesic, Comput. Aided Des. 36, (2004), 447-459.

[27] Williams M.Z., Stein F.M., A triple product of vectors in four-space, Math. Mag., 37, (1964), 230–235.

[28] Willmore T.J., An Introduction to Differential Geometry, Clarendon Press, Oxford, 1959.

[a,b] Ondokuz Mayıs University, Faculty of Arts and Sciences, Department of Mathematics, 55139 Kurupelit/Samsun TURKEY

E-mail address: [a] erginbayram@yahoo.com

E-mail address: [b] kasape@omu.edu.tr